\documentclass{amsart}
\usepackage{palatino,hhline, bbm}
\usepackage{arydshln} 
\usepackage{amsthm, amsmath}
\usepackage{amssymb} 
\usepackage{amscd}
\usepackage{mathrsfs}
\renewcommand{\mathcal}{\mathscr}
\usepackage[hmargin=3cm, vmargin=3cm]{geometry}
\usepackage[breaklinks,bookmarksopen,bookmarksnumbered]{hyperref}
\usepackage[all]{xy}

\theoremstyle{plain}
\newtheorem{thm}{Theorem}

\newtheorem{prop}{Proposition}

\theoremstyle{remark}

\newcommand\pr{\noindent\textit{Proof} : }

\newcommand\rond{\kern 1pt{\scriptstyle\circ}\kern 1pt}

\newcommand\Pic{\operatorname{Pic}}

\newcommand{\mo}{\smallsetminus}

\newcommand\Z{\mathbb{Z}}

\newcommand\C{\mathbb{C}}
\renewcommand\P{\mathbb{P}}

\newcommand\G{\mathbb{G}}

\def\qfl#1{\buildrel {#1}\over {\longrightarrow}}

\newcommand\iso{\vbox{\hbox to .8cm{\hfill{$\scriptstyle\sim$}\hfill}
\nointerlineskip\hbox to .8cm{{\hfill$\longrightarrow $\hfill}} }}
\newcommand\bir{\vbox{\hbox to .8cm{\hfill{$\scriptstyle\sim$}\hfill}
\nointerlineskip\hbox to .8cm{{\hfill$\dasharrow $\hfill}} }}

\begin{document}
\title{A very general sextic double solid  is not  stably rational}
\author[Arnaud Beauville]{Arnaud Beauville}
\address{Laboratoire J.-A. Dieudonn\'e\\
UMR 7351 du CNRS\\
Universit\'e de Nice\\
Parc Valrose\\
F-06108 Nice cedex 2, France}
\email{arnaud.beauville@unice.fr}

\begin{abstract}
We prove that a double covering of $\P^3$  branched along a very general sextic surface is not stably rational.
\end{abstract}

\maketitle 
\section*{Introduction}
A  projective variety $X$ is \emph{stably rational} if $X\times \P^m$ is rational for some integer $m$.  The paper \cite{V} of C. Voisin introduces a new approach to show that some complex varieties are not stably rational. This is  applied in \cite{V} (resp.\ \cite{B}) to prove that a double covering of $\P^3$ (resp.\ $\P^4$ or $\P^5$) branched along a very general quartic hypersurface is not stably rational, and in \cite{CP} to prove the same result for a very general quartic threefold. 

Using the same approach we will prove:
\begin{thm}
A double covering of $\P^3$  branched along a very general sextic surface is not stably rational.
\end{thm}
These ``sextic double solids" are the Fano  threefolds with Picard number 1 of minimal degree ($-K^3=2$).  They are already known to be non-rational \cite{I}; whether they are unirational or not is  unknown. 

\medskip	
We will use Voisin's method in the following form (\cite{V}, Theorem 1.1 and Remark 1.3):
\begin{prop}\label{def}
Let $B$ be a smooth complex variety, $\mathrm{o}$ a point of $B$, $f:\mathcal{X}\rightarrow B$  a flat, projective morphism, such that  the generic fiber of $f$ is smooth, and that the only singularities of the fiber $X:=\mathcal{X}_{\mathrm{o}}$ are ordinary double points. Assume that for a  desingularization $\tilde{X} $ of $X$,  
 the torsion subgroup of $H^3(\tilde{X},\Z )$ is non trivial.
Then for a very general point $b\in B$, the fiber $\mathcal{X}_b$ is not stably rational.
\end{prop}
Thus to  prove the theorem it suffices to find a nodal sextic surface $\Delta \subset \P^3$ such that the desingularization $\tilde{X} $ of the double cover  $X$ of $\P^3$  branched along $\Delta $  satisfies $\mathrm{Tors}\,H^3(\tilde{X},\Z )\neq 0$. Such a surface is described in \cite{IK}. We give here another construction, perhaps simpler; it is not clear to us how the two constructions are related.

As in \cite{B} (and in \cite{IK}) we use a family of quadric surfaces over $\P^3$, with discriminant locus $\Delta $ of degree 6; the quadric fibration provides a natural $\P^1$-bundle over $X_{sm}$, and this gives a 2-torsion class in $H^3(X_{sm},\Z)$, which extends to $H^3(\tilde{X},\Z )$. To construct our quadric fibration we 
start from a cubic fivefold  $V\subset\P^5$, and project from a 2-plane contained in $V$. We show that the associated $\P^1$-bundle has no rational section (Proposition \ref{sec}), and that this provides a nonzero 2-torsion class in  $H^3(\tilde{X} ,\Z)$ (Proposition \ref{H3}). 

\medskip	
{\small I am  indebted to C. Shramov for pointing out the paper \cite{IK}, and to A. Collino for spotting an inaccuracy in the first version of this note.}

\bigskip	
\section*{The construction}
We work over $\C$. Let $V\subset \P^6$ be a smooth cubic fivefold, and   $P\subset V$ a 2-plane. We choose coordinates $(X_0,\ldots ,X_2;\allowbreak Y_0,\ldots ,Y_3)$ on $\P^6$ such that $P$ is given by $Y_0=\ldots =Y_3=0$ and $V$ by
\begin{equation}
\sum_{ i,j}A_{ij}X_iX_j +\sum_iB_iX_i+C=0\label{q}\end{equation} 
where  $A_{ij}, B_i,C$ are homogeneous forms in $(Y_0,\ldots ,Y_3)$ of degree 1, 2 and 3.

Let $\hat{V}$ denote the variety obtained by blowing up $V$ along $P$. The projection from $P$ defines a rational map $V\dasharrow \P^3$, which extends to a morphism $q:\hat{V}\rightarrow \P^3$. The fiber of $q$ at a point $y=(Y_0,\ldots ,Y_3)$ of $\P^3$ is the projective completion of the quadric in $\mathbb{A}^3$ defined by equation (\ref{q}). 

Let $\Delta \subset \P^3$ be the discriminant surface of the quadric fibration $q$, that is, the locus of points $y\in\P^3$ such that $q^{-1}(y)$ is singular. It is defined by the 6th degree equation
\[
\det\left(
\begin{array}{c:c}
\raisebox{15pt}{{\Large\mbox{ {$(A_{ij})$} }}} &\raisebox{15pt}{\mbox{$(B_i)$}}\\[-2ex]
 \hdashline\\[-2.5ex]
 (B_i) & C
\end{array}
\right)=0\ .
\]

According to \cite{C}, Thm. 2.2, for a general choice of the forms $A_{ij}, B_i$ and $C$, the surface $\Delta $ is smooth except for a finite set $\Sigma \subset \Delta $ of $31$ ordinary double points. We will assume from now on that this condition holds. 
The quadric $q^{-1}(y)$ has rank 3 for $y\in \Delta \mo\Sigma $, and rank 2 for $y\in\Sigma $.  

Let $\pi :X\rightarrow \P^3$
 be the double covering branched along $\Delta $. Then $X$ is smooth except for the $31$ ordinary double points lying above $\Sigma $. 
 The generatrices of the quadric $q^{-1}(y)$ 
 are parametrized by two disjoint rational curves for $y\in \P^3\mo\Delta  $, one rational curve for $y\in \Delta \mo\Sigma$. This defines a $\P^1$-bundle $\varphi : G\rightarrow X_{sm} $ onto the smooth locus of $X$. 
 
\begin{prop}\label{sec}
$a)$ The fibration $q:\hat{V}\rightarrow \P^3$ admits no rational section.

$b)$ The $\P^1$-bundle $\varphi $ admits no rational section.
\end{prop}
\pr $a)$ If $q$ admits a rational section,  the closure of its image is a subvariety $Z$ of $\hat{V}$ whose class $[Z]\in H^4(\hat{V},\Z)$ satisfies  $([Z]\cdot q^*y)=1$ for $y\in\P^3$. Let us show that this is impossible.

Consider the blowing-up 
\[\xymatrix@M=6pt{E\ar[d]^p \ar@{^(->}[r]^i & \hat{V}\ar[d]^b  &\\
P\ar@{^(->}[r] & V&\hskip-2.3cm .
}\]
The exceptional divisor $E$ is the hypersurface in $P\times \P^3$ given by $\sum A_{ij}(y)X_iX_j=0$; the projections of $E$ onto $P$ and $\P^3$ are $p$ and $q':=q\rond i$.
The group $H^2(E,\Z)$ is generated by the classes $p^*\ell$ and $q'^*\pi $, where $\ell$  is the class of a line in $H^2(P,\Z)$ and $\pi $ the class of a plane in $\P^3$. 
Let $h\in H^2(V,\Z)$ be the  class of a hyperplane section of $V$; the group $H^4(\hat{V},\Z)$ is generated by the classes $b^*h^2$, $i_*p^*\ell$ and $i_*q'^*\pi $ (see e.g. \cite{B0}, Prop. 0.1.3).  

Let   us compute the intersection number of these classes with the fiber of $q$ at a point $y\in\P^3$. The class $b^*h^2$  induces on the quadric $q^{-1}(y)$ the intersection with a line, hence $(b^*h^2\cdot q^*y)=2$.
For $d\in H^2(E,\Z)$, we have $(i_*d\cdot q^*y)=(d\cdot q'^*y)$. This is zero for $d=q'^*\pi $. The class $p^*\ell$  is the class of a line $\sum a_iX_i=0$, so its  intersection  with the conic $q'^{-1}(y)$ consists of two points, hence $(i_*p^*\ell\cdot q^*y)=2$. It follows that $(\alpha \cdot q^*y)$ is even for any $\alpha \in H^4(\hat{V},\Z)$, so $q$ does not admit a rational section.

\medskip	
$b)$ Suppose $\varphi $ admits a section over some Zariski open subset of $X_{sm}$; we can assume that this subset is of the form $\pi ^{-1}(W)$ for some Zariski open subset $W$ of $\P^3\mo \Delta $. For $w\in W$, the section  maps the two points of $\pi ^{-1}(w)$ to 
two generatrices  of the quadric $q^{-1}(w)$, one in each system.  These two generatrices intersect in one point $s(w) $ of the quadric. This gives a rational section $s$ of $q$, thus contradicting $a)$.\qed

\medskip	

Let $\tilde{X} \rightarrow X$ be the resolution obtained by blowing up the double points; the exceptional divisor $Q$ is a disjoint union of 31 smooth quadrics. 
\begin{prop}\label{H3}
The $2$-torsion subgroup of $H^3(\tilde{X},\Z )$ is nonzero.
\end{prop}
\pr Put $U:=\tilde{X}\mo Q\cong X_{sm} $. The Gysin exact sequence
\[H^1(Q,\Z)\rightarrow H^3(\tilde{X},\Z )\rightarrow H^3(U,\Z)\rightarrow H^2(Q,\Z)\]
shows that the restriction map induces an isomorphism on the torsion subgroups of $H^3(-,\Z)$. Thus it suffices to prove the statement for $H^3(U,\Z )$.

The $\P^1$-bundle $\varphi $ gives a  class $[\varphi ]$ in the 2-torsion subgroup  $\mathrm{Br}_2(U)$ of the Brauer group of $U$; the assertion $b)$ of Proposition \ref{sec} means that this class is nonzero.  
Let us recall how such a class gives a 2-torsion class in $H^3(U,\Z)$, the topogical Brauer class (see \cite{B}, or \cite{N}, 1.1). 
The exact sequences   $\ 0\rightarrow \{\pm 1\}  \rightarrow \G_m\rightarrow \G_m\rightarrow 0\ $ (for the \'etale topology) and $0\rightarrow \Z \qfl{\times 2}\Z\rightarrow \Z/2\rightarrow 0 $ (for the classical topology) give rise to a commutative diagram of exact sequences
\[\xymatrix{\Pic(U)\ar[r]\ar[d]_{c_1}  & H^2(U,\Z/2)\ar[r]\ar@^{=}[d]  & \mathrm{Br}_2(U)  \longrightarrow  0 \\
H^2(U,\Z)\ar[r] & H^2(U,\Z/2)\ar[r]^{\partial} & H^3(U,\Z)\ .\phantom{\longrightarrow}}\]

Therefore $\partial $ induces a homomorphism $\bar{\partial }:\mathrm{Br}_2(U)\rightarrow H^3(U,\Z)$,
which is injective if $c_1:\Pic(U)\rightarrow \allowbreak H^2(U,\Z)$ is surjective. This is indeed the case:   in the commutative diagram
\[\xymatrix{ \Pic(\tilde{X} )\ar[r]^{c_1}\ar[d] & H^2(\tilde{X},\Z )\ar[d]\ \\
 \Pic(U) \ar[r]^{c_1} & H^2(U,\Z )\ , }\] the top horizontal arrow is surjective because $H^2(\tilde{X} ,\mathcal{O}_{\tilde{X} })=0$;  the restriction map $H^2(\tilde{X},\Z )\rightarrow \allowbreak H^2(U,\Z )$ is surjective because of the Gysin exact sequence
$\ H^2(\tilde{X} ,\Z)\rightarrow  H^2(U,\Z)\rightarrow H^1(Q,\Z) =0$. Thus  $\bar{\partial }([\varphi ])$ is a nonzero 2-torsion class in $H^3(U,\Z)$, hence the Proposition.\qed

\medskip	
Theorem 1 follows by   taking for $B$ the space of sextic surfaces in $\P^3$, for $\mathcal{X}$ the family of double coverings of $\P^n$ branched along those surfaces, and
for $\mathrm{o}\in B$ the point corresponding to the discriminant surface $\Delta $.

\end{document}